\documentstyle[11pt,leqno,a4]{article}
\begin{document}
$$$$
\vspace{2cm}
\begin{center}
{\Large \bf On continuous solutions of a problem of R.Schilling}
\end{center}
\vspace{1ex}
\begin{center}
{\large \bf Janusz Morawiec}
\end{center}
\normalsize
\begin{center}
Dedicated to Professor J\'{a}nos Acz\'{e}l with best wishes on
his 70th birthday
\end{center}
\vspace{1ex}
\begin{quotation}
\center
\footnotesize
\begin{flushleft}
{\bf Abstract.} It is proved that if
\end{flushleft}
$q\in\{(\sqrt{3}-1)/2,(3-\sqrt{5})/2,\sqrt{2}-1,(\sqrt{5}-1)/2\}$
\vspace{-1.5ex}
\begin{flushleft}
then\hfill the\hfill zero\hfill function\hfill is\hfill the\hfill only
\hfill solution\hfill $f\! :\! {\rm I\! R}\rightarrow {\rm I\! R}$
\hfill of\hfill (1)\hfill satisfying\hfill (2)\hfill and\hfill
\hfill right-hand-side\hfill or\\
left-hand-side\hfill continuous\hfill at\hfill each\hfill point
\hfill of\hfill the\hfill interval\hfill $(-q/(1-q),-q/(1-q)+\delta )$
\hfill or\hfill of\hfill the\hfill interval\\
$(q/(1-q)-\delta ,q/(1-q))$ with some $\delta >0$.
\end{flushleft}
\end{quotation}
\normalsize
\vspace{1ex}

Studies of a physical problem have led R.Schilling to the
functional equation
\begin{equation}
f(qx)=\frac{1}{4q}[f(x-1)+f(x+1)+2f(x)],
\end{equation}
where $q$ is a fixed number from the open interval $(0,1)$, and
to its solutions $f\! :\! {\rm I\! R}\rightarrow {\rm I\! R}$
such that
\begin{equation}
f(x)=0 \hspace{5ex} {\rm for} \hspace{3ex} \mid x \mid >\frac{q}{1-q}.
\end{equation}

In what follows any solution $f\! :\!{\rm I\! R}\rightarrow {\rm I\! R}$
of (1) satisfying (2) will be called a {\it solution of Schilling's
problem}.

Except the cases $q\in\{ 2^{-\frac{1}{n}}:n\in{\rm I\!N}\}$ we do not
know any nonzero solution of Schilling's problem.
(If $q=2^{-\frac{1}{n}}$ then the functions
$f_{1}(x)=\max \{0,1-|x|\}$ and
$f_{n}=f_{1}(x)*f_{1}(2^{-\frac{1}{n}}x)*\cdots*
f_{1}(2^{-\frac{n-1}{n}}x)$,
respectively, are continuous solutions of this problem.) A partial
explanation of this contains the following theorem of K.Baron and
P.Volkmann [3]: {\it The linear space of the Lebesgue integrable
solutions of Schilling's problem is at most one-dimensional.}
The same cencerns Riemann integrable solutions of Schilling's problem
(see [4; Corollary 3] by W.F\"{o}rg-Rob). The reader interested in
further results up-to-now obtained is referred to the wide paper [4] by
W.F\"{o}rg-Rob, [1] by K.Baron, [2] by K.Baron, A.Simon and
P.Volkmann and to [5] and [6].

It follows from [2; Th\'{e}or\`{e}me 4] by K.Baron, A.Simon and
P.Volkmann that for $q\in (0,1/2)$ there are no nonzero continuous
solutions of Schilling's problem.  The same concerns the case
$q=(\sqrt{5}-1)/2$ and the solutions which are integrable and continuous
at zero (see [2 ; Th\'{e}or\`{e}me 5]).

In the present paper we are interested in solutions which are
right-hand-side
or left-hand-side continuous at each point of one of the intervals:
$$
(-Q,-Q+\delta),\hspace{4ex}(Q-\delta,Q),
$$
where
$$
Q=\frac{q}{1-q}
$$
and $\delta$ is a positive real number, in the case where $q$ is one of
the numbers:
\begin{equation}
\frac{\sqrt{3}-1}{2},\hspace{4ex}\frac{3-\sqrt{5}}{2},\hspace{4ex}
\sqrt{2}-1,\hspace{4ex}\frac{\sqrt{5}-1}{2}.
\end{equation}

We start with the following lemma which follows from [4; Lemmas 1(ii)
and 2(ii)].

{\bf Lemma}. {\it If a solutions of Schilling's problem vanishes either
on an interval $(Q-\delta , Q]$ or on an interval $[-Q,-Q+\delta )$ for
some $\delta >0$, then it vanishes everywhere.}

Our main result reads as follows.

{\bf Theorem}. {\it If $q$ is one of the numbers {\rm (3)}, then any
solution of Schilling's problem vanishes on the set
${\rm Z\!\!Z}+q{\rm Z\!\!Z}$.}

{\bf Proof}. Let
$$
q=\frac{\sqrt{3}-1}{2}.
$$
Then
\begin{equation}
1-2q-2q^{2}=0.
\end{equation}

Denote by $A_{0}$ the set of all the numbers of the form
$$
\varepsilon (m+2nq),\hspace{4ex}\varepsilon (m+1+n-2nq),
$$
where $\varepsilon\in\{-1,1\}$ and $m$, $n$ are non-negative integers,
and put
$$
A_{n}=\{\varepsilon (m-2nq):\hspace{1ex}
m\in\{1,\cdots ,n\},\varepsilon\in\{-1,1\}\}\hspace{5ex}{\rm for}
\hspace{3ex} n\in {\rm I\! N}.
$$
Then
\begin{equation}
\bigcup_{n=0}^{\infty}A_{n}={\rm Z\!\!Z}+2q{\rm Z\!\!Z}.
\end{equation}
Using induction we shall show that
\begin{equation}
f(x)=0\hspace{5ex}{\rm for}\hspace{3ex} x\in A_{n}
\end{equation}
and for every non-negative integer $n$.

To get (6) for $n=0$ let us observe (cf. also [5; Remark 2(ii)]) that
$q<1/2$ gives
\begin{equation}
f(0)=0.
\end{equation}
Further, if at least one of the non-negative integer $m$, $n$ is
positive then $m+2nq\geq 2q>Q$ which jointly with (2) gives
$f(\varepsilon (m+2nq))=0$ for $\varepsilon\in\{-1,1\}$. Moreover, if
$m$ and $n$ are positive integers then $m+n-2nq\geq 1+(1-2q)>1>Q$ whence
$f(\varepsilon (m+n-2nq))=0$ for $\varepsilon\in\{-1,1\}$. This proves
that $f$ vanishes on $A_{0}$.

Fix  now a positive integer $n$ and assume that $f$ vanishes on
$A_{0},\cdots A_{n-1}$. Let $x_{0}\in A_{n}$.
Then $x_{0}=\varepsilon (m-2nq)$ for some $\varepsilon\in \{-1,1\}$ and
$m\in\{1,\cdots ,n\}$. Hence and from (4) we get
$$
x_{0}=\varepsilon [(2q+2q^{2})m-2nq]=2q\varepsilon [(m-n)+mq].
$$
Putting now $x=x_{0}/q$ into (1) we have
\begin{equation}
\hspace{-42ex}f(x_{0})=\frac{1}{4q}[f(x-1)+f(x+1)+2f(x)]
\end{equation}
$$
\hspace{15ex}=\frac{1}{4q}[f(\varepsilon [2(m-n)-1+2mq])+
f(\varepsilon [2(m-n)+1+2mq])+2f(\varepsilon [2(m-n)+2mq])].
$$
Moreover, each of the points
$$
\varepsilon [2(m-n)-1+2mq],\hspace{4ex}
\varepsilon [2(m-n)+1+2mq],\hspace{4ex}
\varepsilon [2(m-n)+2mq]
$$
belongs to one of the sets $A_{0}$, $A_{m}$.
Hence, if $n\geq 2$ and $m<n$, then using (8) and the induction
hypothesis we obtain $f(x_{0})=0$. In particular,
\begin{equation}
f(\varepsilon (1-2nq))=0.
\end{equation}
If $n\geq 2$ and $m=n$, then (8), (2) and (9) give
$$
f(x_{0})=\frac{1}{4q}f(\varepsilon (-1+2nq))=0.
$$
Finally assume that $n=m=1$. Then, according to (8) and (2),
$$
f(x_{0})=\frac{1}{4q}f(\varepsilon (-1+2q)),
$$
i.e.,
$$
f(\varepsilon (1-2q))=\frac{1}{4q}f(\varepsilon (-1+2q)).
$$
Consequently, since the above equality holds for every
$\varepsilon\in\{-1,1\}$,
$$
f(\varepsilon (1-2q))=\frac{1}{(4q)^{2}}f(\varepsilon (1-2q))
$$
and so $f(x_{0})=0$. This ends the induction proof of (6) and jointly
with (5) shows that
\begin{equation}
f\hspace{1ex}{\rm vanishes\hspace{1ex} on}\hspace{1ex}
{\rm Z\!\!Z}+2q{\rm Z\!\!Z}.
\end{equation}

To prove that $f$ vanishes on ${\rm Z\!\!Z}+q{\rm Z\!\!Z}$ let us
observe that taking into account (4) we have
$$
{\rm Z\!\!Z}+q{\rm Z\!\!Z}=(2q+2q^{2}){\rm Z\!\!Z}+q{\rm Z\!\!Z}
\subset q({\rm Z\!\!Z}+2q{\rm Z\!\!Z})
$$
which together with (1) and (10) ends the proof in the case where $q$
is the first of the numbers (3).

Now we pass to the case
$$
q=\frac{3-\sqrt{5}}{2}.
$$
Observe that
\begin{equation}
1-3q+q^{2}=0\hspace{4ex}{\rm and}\hspace{4ex}Q=1-q.
\end{equation}

Denote now by $A_{1}$ the set of all the numbers of the form
$$
\varepsilon (m+nq),\hspace{4ex}\varepsilon (m+1+n-(n+1)q),
$$
where $\varepsilon\in\{-1,1\}$ and $m$, $n$ are non-negative integers,
and put
$$
A_{n}=\{\varepsilon (m-nq):\hspace{1ex}
m\in\{1,\cdots ,n-1\},\varepsilon\in\{-1,1\}\}\hspace{5ex}{\rm for}
\hspace{3ex} n\geq 2.
$$
Then
\begin{equation}
\bigcup_{n=1}^{\infty}A_{n}={\rm Z\!\!Z}+q{\rm Z\!\!Z}.
\end{equation}
We shall show that (6) holds for every positive integer $n$.

To obtain (6) for $n=1$ let us notice that $1/4\neq q<1/2$ which jointly
with the second part of (11) (cf. [5; Remarks 1 and 2]) gives
\begin{equation}
f(0)=f(1-q)=f(-1+q)=0.
\end{equation}
Moreover, if at least one of the non-negative integer $m$, $n$ is
positive then $m+n+1-(n+1)q\geq 1>Q$. Thus, according to (2),
$f(\varepsilon (m+n+1-(n+1)q))=0$ for $\varepsilon\in\{-1,1\}$.
Now observe that putting $x=\varepsilon\in\{ -1,1\}$ into (1) and
applying (2) and (13) we get $f(\varepsilon q)=0$. Further, if $m\geq 1$
or $n\geq 2$, then $m+nq\geq 2q>Q$ which jointly with (2) shows that
$f(\varepsilon (m+nq))=0$ for $\varepsilon\in\{-1,1\}$. This proves
that $f$ vanishes on $A_{1}$.

Fix  now an integer $n\geq 2$ and assume that $f$ vanishes on
$A_{1},\cdots A_{n-1}$. Let $x_{0}\in A_{n}$.
Then $x_{0}=\varepsilon (m-nq)$ for some $\varepsilon\in \{-1,1\}$ and
$m\in\{1,\cdots ,n-1\}$. Hence and from (11) we get
$$
x_{0}=\varepsilon [(3q-q^{2})m-nq]=q\varepsilon [(3m-n)-mq].
$$
Putting now $x=x_{0}/q$ into (1), using the induction hypothesis and the
inclusion $\varepsilon ({\rm Z\!\!Z}-mq)\subset A_{m}\cup A_{1}$ we have
$$
\hspace{-38ex}f(x_{0})=\frac{1}{4q}[f(x-1)+f(x+1)+2f(x)]
$$
$$
\hspace{15ex}=\frac{1}{4q}[f(\varepsilon (3m-n-1-mq))+
f(\varepsilon (3m-n+1-mq))+2f(\varepsilon (3m-n-mq))].
$$
This ends the induction proof of (6) and jointly with (12) shows that
$f$ vanishes on ${\rm Z\!\!Z}+q{\rm Z\!\!Z}$.

Let us pass to the case
$$
q=\sqrt{2}-1
$$
and notice that now
\begin{equation}
1-2q-q^{2}=0.
\end{equation}
Moreover, since $q<1/2$, (7) holds. Consider the sets $A_{n}$ introduced
in the previous case. Now putting in (1) $x=\varepsilon\in\{ -1,1\}$
and using (2) and (7) we obtain
\begin{equation}
f(\varepsilon q)=0,
\end{equation}
whereas putting in (1) $x=\varepsilon (1+q)$ and using (2), (15) and
(14) we have $f(\varepsilon (1-q))=0$. Moreover, if $m\geq 1$ or
$n\geq 2$, then $m+nq\geq 2q>Q$ which jointly with (2) gives
$f(\varepsilon (m+nq))=0$ for $\varepsilon\in\{-1,1\}$. Besides,
if $m\geq 1$ or $n\geq 1$, then $m+n+1-(n+1)q\geq 2-2q>1>Q$, and
applying (2) we obtain $f(\varepsilon (m+n+1-(n+1)q))=0$
for $\varepsilon\in\{-1,1\}$. Thus $f$ vanishes on $A_{1}$. Now using
induction and (14) (as we do it in the case above) we get that $f$
vanishes on the set (12).

Finally we assume that
$$
q=\frac{\sqrt{5}-1}{2}.
$$
Observe that now
\begin{equation}
1-q-q^{2}=0\hspace{4ex}{\rm and}\hspace{4ex}Q=\frac{1}{q}=1+q.
\end{equation}

Denote by $A_{1}$ the set of all the numbers of the form
$$
\varepsilon (m+nq),\hspace{4ex}\varepsilon (m+1+n-nq),
\hspace{4ex}\varepsilon (1-q),
$$
where $\varepsilon\in\{-1,1\}$ and $m$, $n$ are non-negative integers,
and put
$$
A_{n}=\{\varepsilon (m-nq):\hspace{1ex}
m\in\{1,\cdots ,n\},\varepsilon\in\{-1,1\}\}\hspace{5ex}{\rm for}
\hspace{3ex} n\geq 2.
$$
Evidently (12) holds.

Now we shall show that (6) is fulfiled for $n=1$. To this end notice
first that $q\neq 1/4$, whence according to [5; Remarks 1 and 2(i)] and
(16) we have
\begin{equation}
f(1+q)=f(-1-q)=0.
\end{equation}
Moreover,
\begin{center}
if $m\geq 1$ and $n\geq 1$, then $m+nq\geq 1+q=Q$,\\
if $n\geq 3$, then $nq\geq 3q>Q$ and $n-q\geq 3-q>Q$,
\end{center}
and
\begin{center}
if $m\geq 0$ and $n\geq 2$, then $m+1+n-nq\geq 3-2q>Q$.
\end{center}
This jointly with (2) gives that $f$ vanishes on the set
$$
A_{1}\setminus \{ 0,1,-1,q,-q,2q,-2q,1-q,-1+q,2-q,-2+q\}.
$$
Putting in (1) in turn $x=0$, $x=\varepsilon $, $x=\varepsilon /q$ with
$\varepsilon\in \{-1,1\}$ and making use of (2), (17) and(16) we have
\begin{equation}
(4q-2)f(0)=f(\varepsilon )+f(-\varepsilon ),
\end{equation}
\begin{equation}
4qf(\varepsilon q)=2f(\varepsilon )+f(0),
\end{equation}
\begin{equation}
4qf(\varepsilon )=f(\varepsilon q).
\end{equation}
Taking into account (19) and (20) we get
$$
f(\varepsilon )=f(-\varepsilon )
$$
and using (18), (19) and (20) again we obtain (after some calculations)
\begin{equation}
f(0)=f(1)=f(-1)=f(q)=f(-q)=0.
\end{equation}

Putting in (1): $x=2\varepsilon $ and $x=\varepsilon (2/q-1)$ with
$\varepsilon\in \{-1,1\}$ and using (21), (2) and (16) we have
\begin{equation}
f(2q)=f(-2q)=f(2-q)=f(-2+q)=0.
\end{equation}
Finally putting $x=\varepsilon q$ into (1) and making use of (16), (17)
and (21) we see that
$$
f(1-q)=f(-1+q)=0
$$
which jointly with (21) and (22) gives that $f$ vanishes on the whole
of the set $A_{1}$.

Fix now an $n\geq 2$ and assume that $f$ vanishes on
$A_{1},\cdots ,A_{n-1}$. Let $x_{0}\in A_{n}$.
Then $x_{0}=\varepsilon (m-nq)$ for some $\varepsilon\in \{-1,1\}$ and
$m\in\{1,\cdots ,n\}$. Hence and from (16) we get
$$
x_{0}=\varepsilon [(q+q^{2})m-nq]=q\varepsilon [(m-n)+mq].
$$
Putting now $x=x_{0}/q$ into (1) we obtain
\begin{equation}
\hspace{-37ex}f(x_{0})=\frac{1}{4q}[f(x-1)+f(x+1)+2f(x)]
\end{equation}
$$
\hspace{15ex}=\frac{1}{4q}[f(\varepsilon [(m-n)-1+mq])+
f(\varepsilon [(m-n)+1+mq])+2f(\varepsilon [(m-n)+mq])]
$$
and each of the points
$$
\varepsilon [(m-n)-1+mq],\hspace{4ex}
\varepsilon [(m-n)+1+mq],\hspace{4ex}
\varepsilon [(m-n)+mq]
$$
belongs to one of the sets $A_{1}$, $A_{m}$. Hence, if $m<n$, then on
account of the induction hypothesis we have $f(x_{0})=0$. In particular,
$f(\varepsilon (1-nq))=0$. If $m=n$, then using (23), (2), (22) and the
above fact we see that
$$
f(x_{0})=\frac{1}{4q}[f(\varepsilon (-1+nq))+2f(\varepsilon nq)]=0
$$
and the proof is completed.

{\bf Corollary.} {\it If $q$ is one of the numbers {\rm (3)} then the
zero function is the only solution of Schilling's problem which is
either right-hand-side or left-hand-side continuous at each point of
one of the intervals $(-Q,-Q+\delta )$, $(Q-\delta ,Q)$ where
$\delta $ is a positive real number.}
\vspace{3ex}

{\bf Acknowledgement.} This research was supported by the State
Committee for Scientific Research Grant No. 2 1062 91 01.

\vspace{3ex}
\begin{flushleft}
{\bf References}
\end{flushleft}
\footnotesize
[1] K. Baron, {\it On a problem of R.Schilling.} Berichte der
Mathematisch-statistischen Sektion in der Forschungsge-\linebreak
\mbox{\hspace{3ex}} sellschaft Joanneum--Graz, Bericht Nr. 286
(1988).\newline
[2] K. Baron, A.Simon et P.Volkmann, {\it Solutions d'une
\'{e}quation functionnelle dans l'espace des distributions
temp\'{e}r\'{e}es},\linebreak
\mbox{\hspace{3ex}} Manuscript.\newline
[3] K. Baron et P.Volkmann, {\it Unicit\'{e} pour une \'{e}quation
functionnelle}, Rocznik Naukowo-Dydaktyczny WSP w\linebreak
\mbox{\hspace{3ex}} Krakowie. Prace Matematyczne 13 (1993),
53-56.\newline
[4] W.F\"{o}rg-Rob, {\it On a problem of R.Schilling}, Mathematica
Pannonica, to appear.\newline
[5] J.Morawiec, {\it On bounded solutions of a problem of R.Schilling},
Annales Mathematicae Silesianae, to appear.\newline
[6] J.Morawiec, {\it On a problem of R.Schilling}, Report of
Meeting: The Thirty-first International Symposium on Func-\
\mbox{\hspace{3ex}} tional Equations (August 22 - August 28, 1993,
Debrecrn, Hungary), Aequationes Mathematica 47(1994), 284-285.
\vspace{3ex}
\normalsize

\begin{flushleft}
Instytut Matematyki,\\
Uniwersytet \'{S}l\c aski,\\
ul. Bankowa 14,\\
PL-40-007 Katowice
\end{flushleft}
\begin{flushright}
Eingegangen am 24. Juni 1994
\end{flushright}
\end{document}